\documentclass[12pt]{amsart}
\usepackage{amsfonts, amsmath, amssymb, latexsym}
\usepackage[all]{xy}
    \title[Co-$t$-structures]{Compact corigid objects in triangulated Categories and co-$t$-structures}
    \author{David Pauksztello}
    \date{14th October 2007}

\newtheorem{dfn}{Definition}[section]

\newtheorem{thrm}[dfn]{Theorem}

\newtheorem{lem}[dfn]{Lemma}

\newtheorem{prp}[dfn]{Proposition}

\newtheorem{corr}[dfn]{Corollary}

\newtheorem{rmk}[dfn]{Remark}

\newtheorem{set}[dfn]{Setup}

\newtheorem{egg}[dfn]{Example}

\newenvironment{prf}{\noindent\textbf{Proof: }}{$\Box$\ }

\newcommand{\comment}[1]{}

\newcommand{\D}{\mathcal{D}}

\newcommand{\T}{\mathcal{T}}

\newcommand{\A}{\mathcal{A}}

\newcommand{\F}{\mathcal{F}}

\newcommand{\Hom}[3]{\textnormal{Hom}_{#1}(#2,#3)}

\newcommand{\rightiso}{\stackrel{\sim}{\longrightarrow}}

\newcommand{\hocolim}[1]{\textnormal{hocolim}(#1)}

\newcommand{\Endo}[1]{\textnormal{End}(#1)}

\newcommand{\rightlabel}[1]{\stackrel{#1}{\longrightarrow}}

\newcommand{\Z}{\mathbb{Z}}

\newcommand{\df}[1]{\begin{dfn}\emph{#1}}

\newcommand{\edf}{\end{dfn}}

\newcommand{\thm}{\begin{thrm}}

\newcommand{\ethm}{\end{thrm}}

\newcommand{\lemma}{\begin{lem}}

\newcommand{\elemma}{\end{lem}}

\newcommand{\prop}{\begin{prp}}

\newcommand{\eprop}{\end{prp}}

\newcommand{\rem}[1]{\begin{rmk}\emph{#1}}

\newcommand{\erem}{\end{rmk}}

\newcommand{\pf}{\begin{prf}}

\newcommand{\epf}{\end{prf}}

\newcommand{\cor}{\begin{corr}}

\newcommand{\ecor}{\end{corr}}

\newcommand{\setup}[1]{\begin{set}\emph{#1}}

\newcommand{\esetup}{\end{set}}

\newcommand{\eg}[1]{\begin{egg}\emph{#1}}

\newcommand{\eeg}{\end{egg}}

\newcommand{\tri}[3]{#1\rightarrow #2\rightarrow #3\rightarrow \Sigma #1}

\newcommand{\trilabel}[4]{#1\stackrel{#4}{\longrightarrow} #2\longrightarrow 
#3\longrightarrow \Sigma #1}

\newcommand{\backtri}[3]{\Sigma^{-1}#1\rightarrow #2\rightarrow 
#3\rightarrow #1}

\newcommand{\backtrilabel}[4]{\Sigma^{-1}#1\stackrel{#4}{\longrightarrow} 
#2\longrightarrow #3\longrightarrow #1}

\entrymodifiers={+!!<0pt,\fontdimen22\textfont2>}

\numberwithin{equation}{section}

\begin{document}
	
\address{Department of Pure Mathematics, University of Leeds,
Leeds. LS2 9JT United Kingdom}
\email{davidp@maths.leeds.ac.uk}
\urladdr{http://www.maths.leeds.ac.uk/\~{ }davidp}    
    

    
\begin{abstract}
In the work of Hoshino, Kato and Miyachi, \cite{Miyachi}, the authors look at $t$-structures induced by a compact object, $C$, of a triangulated category, $\T$, which is rigid in the sense of Iyama and Yoshino, \cite{Iyama}. Hoshino, Kato and Miyachi show that such an object yields a non-degenerate $t$-structure on $\T$ whose heart is equivalent to $\textnormal{Mod}(\Endo{C}^{\textnormal{op}})$. Rigid objects in a triangulated category can the thought of as behaving like chain differential graded algebras (DGAs).

Analogously, looking at objects which behave like cochain DGAs naturally gives the dual notion of a corigid object. Here, we see that a compact corigid object, $S$, of a triangulated category, $\T$, induces a structure similar to a $t$-structure which we shall call a co-$t$-structure. We also show that the coheart of this non-degenerate co-$t$-structure is equivalent to $\textnormal{Mod}(\Endo{S}^{\textnormal{op}})$, and hence an abelian subcategory of $\T$.
\end{abstract}

\maketitle

\setcounter{section}{-1}

\section{Introduction}\label{introduction}

Suppose $\T$ is a triangulated category with set indexed coproducts and let $\Sigma:\T\rightarrow \T$ denote its suspension functor. Hoshino, Kato and Miyachi, in \cite{Miyachi}, show that a natural $t$-structure is induced on $\T$ by a suitably  nice compact object of $\T$. In particular, they consider a compact object $S$ of $\T$ which satisfies the following two conditions:

$(1)$ $\Hom{\T}{S}{\Sigma^{i}S}=0 \textnormal{ for all } i>0$;

$(2)$ $\{\Sigma^{i}S \,|\,i\in\mathbb{Z}\}$ is a generating set for $\T$.

\noindent Following the terminology of Iyama and Yoshino, we refer to an object satisfying the first of the two conditions above as \textit{rigid}, see \cite{Iyama}. We shall give precise definitions of the notions of $t$-structure, compact object and generating set in sections \ref{definitions} and \ref{examples}.

If $S$ is a compact rigid object of $\T$ and $\{\Sigma^{i}S \,|\,i\in\mathbb{Z}\}$ is a generating set for $\T$, then the two halves of the $t$-structure obtained in \cite{Miyachi} are given by
\begin{eqnarray*}
\mathcal{X} & = & \{X\in\T \ | \ \Hom{\T}{S}{\Sigma^{i}X}=0 \textnormal{ for } i>0\}, \\
\mathcal{Y} & = & \{X\in\T\ | \ \Hom{\T}{S}{\Sigma^{i}X}=0 \textnormal{ for } i<0\}.
\end{eqnarray*}
This situation bears resemblance to the example of a chain differential graded algebra (DGA) $R$ in its derived category $\D(R)$, whose objects are the differential graded modules over $R$ (DG $R$-modules). Introductions to the theory of DGAs, their derived categories and DG modules can be found in \cite{Aldrich}, \cite{Bernstein} and \cite{FJ2}.

Recall that a DGA $R$ is called a \textit{chain DGA} if $H^{i}(R)=0$ for all $i>0$. Moreover, given a DG $R$-module $M$ we have $$H^{i}(M)\cong \Hom{\D(R)}{R}{\Sigma^{i}M} \textnormal{ for all } i\in\mathbb{Z}.$$ Hence, the object $S$ considered in \cite{Miyachi} is analogous to a chain DGA and the two halves of the $t$-structure it induces are analogous to the full subcategories of DG $R$-modules whose cohomology vanishes in positive and negative degree, respectively.

In the theory of DGAs, when one has a construction for chain DGAs it is natural to ask: what is the dual construction for cochain DGAs? Likewise, it is natural to ask, what is the structure induced by a compact object of a triangulated category which behaves like a cochain DGA?
Recall that a DGA $R$ is called a \textit{cochain DGA} if $H^{i}(R)=0$ for $i<0$. 

Unfortunately, it is well known in the theory of DGAs that constructing a viable dual theory for cochain DGAs is often difficult. In fact, at present the construction of a viable dual theory for DGAs always requires the additional assumption that the DGA $R$ is simply connected in the following sense: $H^{0}(R)$ is a division ring and $H^{1}(R)=0$.
This lack of symmetry between the chain and cochain theories occurs throughout the theory of DGAs and in algebraic topology, see \cite{Avramov}, for example. Thus, we shall consider the case of a compact object of a triangulated category which behaves like a simply connected cochain DGA. 

The structure which is induced by such an object is not a $t$-structure, but it turns out to be almost dual to the notion of a $t$-structure, and as such we call it a co-$t$-structure. Both $t$-structures and co-$t$-structures provide examples of torsion theories in triangulated categories in the sense of Iyama and Yoshino, \cite{Iyama}.
Co-$t$-structures have also been introduced by Bondarko in \cite{Bondarko} where they are called \textit{weight structures}. They are studied in \cite{Bondarko} in connection with the theory of motives and stable homotopy theory.

The paper is organised as follows: in section \ref{definitions}, we recall the concepts of preenvelopes and precovers and set up the notation of perpendicular categories. We then recall the notion of a $t$-structure and introduce the new definition of a co-$t$-structure, about which we prove some elementary properties and compare and contrast this new notion with the existing notion of a $t$-structure. We also say why it is almost dual to a $t$-structure but not exactly dual. In addition, we introduce the definition of the coheart of a co-$t$-structure.

In section \ref{examples}, we look at a canonical example of a co-$t$-structure appearing in the setting of the homotopy category of an additive category. We then look at the simple motivating examples of the $t$-structure induced by a chain DGA and the co-$t$-structure obtained by a simply connected cochain DGA on their respective derived categories. Note again that in order to consider a viable cochain analogue we need to impose the simply connected hypothesis. We present a brief exposition of Hoshino, Kato and Miyachi's theorem, which is obtained in \cite{Miyachi}, which generalises and abstracts the example of the $t$-structure induced on the derived category of a chain DGA. Hoshino, Kato and Miyachi's theorem is presented here as Theorem \ref{miyachitheorem}.

The remainder of the paper is devoted to proving the simply connected cochain analogue of Hoshino, Kato and Miyachi's theorem, which is presented as Theorem \ref{cochaintheorem}, the first half of the proof appearing in sections \ref{preenvelopes} and \ref{particular}. In addition to inducing a $t$-structure, Hoshino, Kato and Miyachi also prove that the heart of the induced $t$-structure, which is well known to be admissible abelian (see \cite{Beilinson} and \cite{Kashiwara}), is equivalent to the module category of the endomorphism algebra of the object inducing the $t$-structure on $\T$. Here, we are able to prove a similar result regarding the coheart of the induced co-$t$-structure. It is known that the coheart of a co-$t$-structure is not always abelian, and may be very rarely so. Indeed, a specific example whose coheart is not abelian is constructed by Bondarko in \cite{Bondarko}. However, the example of a co-$t$-structure which we present in this paper has an abelian coheart by virtue of its equivalence to a module category; this is the subject of section \ref{coheart}.

\section{Definitions, terminology and notation}\label{definitions}

In this section we shall introduce some of the basic definitions, terminology and notation which we shall use throughout this paper. We start by recalling the concept of a preenvelope and the notation of perpendicular categories; then we give the definition of a $t$-structure on a triangulated category and introduce the definition of a co-$t$-structure and some of its basic properties.

Throughout this paper $\T$ will be a triangulated category with set indexed coproducts. We shall denote the suspension functor on $\T$ by $\Sigma:\T\rightarrow \T$, and we shall write $\Hom{}{X}{Y}$ instead of $\Hom{\T}{X}{Y}$ for the Hom-sets of $\T$.

\subsection{Preenvelopes and perpendicular categories}

We recall the following definition from \cite{Enochs}.

\df
{Let $\F$ be a full subcategory of a category $\mathcal{C}$ and suppose $X$ is an object of $\mathcal{C}$. A morphism $\phi:X\rightarrow F$ with $F\in\F$ is called an \textit{$\F$-preenvelope} if for each morphism $X\rightarrow F'$ with $F'\in \F$ there exists a morphism $F\rightarrow F'$ making the following triangle commute.
$$\xymatrix{X\ar[r]^-{\phi}\ar[dr] & F\ar@{-->}[d]^-{\exists} \\
 & F' }$$}
\edf

An $\F$-preenvelope is sometimes called a \textit{left $\F$-approximation}. We obtain the notion of an \textit{$\F$-precover} by dualising the definition above.

\df
{Let $S$ be an object of a triangulated category $\T$. The \textit{subcategory right $n$-perpendicular} to $S$, denoted by $S^{\perp_{n}}$, is given by:}
\begin{equation*}
S^{\perp_{n}}:=\{X\in\T \ | \ \Hom{}{S}{\Sigma^{i}X}=0 \textnormal{ for } i=1,\ldots, n\}.
\end{equation*}
\textnormal{The \textit{subcategory right $\infty$-perpendicular} to $S$, denoted by $S^{\perp_{\infty}}$, is given by:}
\begin{equation*}
S^{\perp_{\infty}}:=\{X\in\T \ | \ \Hom{}{S}{\Sigma^{i}X}=0 \textnormal{ for } i>0\}.
\end{equation*}
\textnormal{Similarly, one can also define the \textit{subcategories left $n$-perpendicular} and \textit{left $\infty$-perpendicular} to $S$.}
\edf

For a subcategory $\mathcal{S}$ of $\T$, we define:
\begin{eqnarray*}
\mathcal{S}^{\perp} & := & \{X\in\T \ | \ \Hom{}{S}{X}=0 \textnormal{ for all } S\in\mathcal{S}\}, \\
{}^{\perp}\mathcal{S} & := & \{X\in\T \ | \ \Hom{}{X}{S} = 0 \textnormal{ for all } S\in\mathcal{S}\}.
\end{eqnarray*}

\subsection{$t$-structures and co-$t$-structures}

The concept of a $t$-structure on a triangulated category $\T$ was first introduced by Beilinson, Bernstein and Deligne in \cite{Beilinson}. The basic theory of $t$-structures can be found in \cite{Beilinson} and \cite{Kashiwara}.

\df
{Let $\T$ be a triangulated category. A pair of full subcategories of $\T$, $(\mathcal{X},\mathcal{Y})$, is called a \textit{$t$-structure} on $\T$ if it satisfies the following properties:}

$(1)$ \textnormal{$\mathcal{X}\subseteq \Sigma^{-1} \mathcal{X}$ and $\Sigma^{-1} \mathcal{Y}\subseteq \mathcal{Y}$;}

$(2)$ \textnormal{$\Hom{}{\mathcal{X}}{\Sigma^{-1}\mathcal{Y}}=0$;}

$(3)$ \textnormal{For any object $Z$ of $\T$ there exists a distinguished triangle $$\tri{X}{Z}{\Sigma^{-1}Y}$$  with $X\in\mathcal{X}$ and $Y\in \mathcal{Y}$.}
\label{t_structure}
\edf

The full subcategories $\mathcal{X}$ and $\mathcal{Y}$ are often denoted by $\T^{t\leqslant 0}$ and $\T^{t\geqslant 0}$, or simply by $\T^{\leqslant 0}$ and $\T^{\geqslant 0}$, respectively; see \cite{Bondarko} and \cite{Miyachi}.

The notion of a $t$-structure has become widespread in the study of triangulated categories and lends itself particularly well to induction arguments in this setting, see for example \cite{Bergh}.

We next introduce the almost dual notion of a co-$t$-structure. We note that co-$t$-structures have also recently been introduced by Bondarko in \cite{Bondarko} where they are called \textit{weight structures}.

\df
{Let $\T$ be a triangulated category. A pair of full subcategories of $\T$, $(\A, \mathcal{B})$, will be called a \textit{co-$t$-structure} on $\T$ if it satisfies the following properties:}

$(0)$ \textnormal{$\A$ and $\mathcal{B}$ are closed under direct summands;}

$(1)$ \textnormal{$\Sigma^{-1}\A\subseteq \A$ and $\mathcal{B}\subseteq \Sigma^{-1}\mathcal{B}$;}

$(2)$ \textnormal{$\Hom{}{\Sigma^{-1}\A}{\mathcal{B}}=0$;}

$(3)$ \textnormal{For any object $X$ of $\T$ there exists a distinguished triangle $$\backtri{A}{X}{B}$$ with $A\in\A$ and $B\in \mathcal{B}$.}
\label{cot}
\edf

In \cite{Bondarko} the full subcategories $\A$ and $\mathcal{B}$ are denoted by $\T^{w\geqslant 0}$ and $\T^{w\leqslant 0}$, respectively.

It is easy to see that $\A$ is closed under direct summands if and only if $\Sigma^{-1}\A$ is closed under direct summands; similarly for $\mathcal{B}$.

One can see that by interchanging the roles of $\mathcal{X}$ and $\mathcal{Y}$ in the definition of a $t$-structure, properties $(1)$, $(2)$ and $(3)$ in Definition \ref{t_structure} become the corresponding properties in the definition of a co-$t$-structure. 
The inclusion of condition $(0)$ in Definition \ref{cot} is the reason why a co-$t$-structure is almost dual to a $t$-structure rather than simply being its dual.
We also note that properties $(2)$ and $(3)$ in Definitions \ref{t_structure} and \ref{cot} make $(\mathcal{X},\Sigma^{-1}\mathcal{Y})$ and $(\Sigma^{-1}\A,\mathcal{B})$ into examples of torsion theories in the sense of \cite{Iyama}.

The notion of a non-degenerate co-$t$-structure can be defined in a manner analogous to that of a non-degenerate $t$-structure.

\df
{A co-$t$-structure $(\A,\mathcal{B})$ on a triangulated category $\T$ will be called \textit{non-degenerate} if we have
$$\bigcap_{n\in\Z}\Sigma^{n}\A=\bigcap_{n\in\Z}\Sigma^{n}\mathcal{B}=\{0\}.$$}
\edf

Recall that a full subcategory $\mathcal{X}$ of a triangulated category $\T$ is said to be \textit{closed under extensions} if, whenever we have a distinguished triangle
$$\tri{X'}{X}{X''}$$
with $X'$ and $X''$ objects of $\mathcal{X}$, then $X$ is also an object of $\mathcal{X}$.  We next give some elementary properties of co-$t$-structures. 

\prop
Let $\T$ be a triangulated category and suppose $(\mathcal{A},\mathcal{B})$ is a co-$t$-structure on $\T$. We have:

(i) For all objects $X$ of $\T$ there exists an $\Sigma^{-1}\A$-precover $\alpha: \Sigma^{-1}A\rightarrow X$.

(ii) For all objects $X$ of $\T$ there exists a $\mathcal{B}$-preenvelope $\beta : X\rightarrow B$.

(iii) We have $\Sigma^{-1}\A = {}^{\perp}\mathcal{B}$ and $\mathcal{B} = (\Sigma^{-1}\A)^{\perp}$.

(iv) $\mathcal{A}$ is closed under extensions.

(v) $\mathcal{B}$ is closed under extensions.
\label{prop3a}
\eprop

\pf
Properties (i) and (ii) are immediate consequences of the definition of a co-$t$-structure: the $\Sigma^{-1}\A$-precover $\alpha: \Sigma^{-1}A\rightarrow X$ and the $\mathcal{B}$-preenvelope $\beta:X\rightarrow B$ are just the first and second morphisms in the distinguished triangle given by property $(3)$ of Definition \ref{cot}. Property (iii) is a consequence of condition $(0)$ and the orthogonality condition $(2)$ of Definition \ref{cot}, and properties (iv) and (v) are easy consequences of property (iii), see \cite{Iyama}.
\epf

\

One sees in Proposition \ref{prop3a} that preenvelopes and precovers replace the truncation functors associated with $t$-structures. In order to obtain the equalities of property (iii), and thus the fact that both halves of the co-$t$-structure are closed under extensions, we need to assume condition $(0)$ of Definition \ref{cot} which says that both halves of a co-$t$-structure are closed under direct summands.

Let $(\mathcal{X},\mathcal{Y})$ be a $t$-structure on a triangulated category $\T$. The intersection $\mathcal{H}=\mathcal{X}\cap\mathcal{Y}$, of both halves of the $t$-structure is called the \textit{heart} of the $t$-structure. It has the nice property that it is an abelian subcategory of $\T$. In particular, the hearts of $t$-structures provide a means of obtaining abelian categories from triangulated categories. We define an analogous notion for co-$t$-structures.

\df
{Let $(\mathcal{A},\mathcal{B})$ be a co-$t$-structure for a triangulated category $\T$. The intersection $\mathcal{C}=\mathcal{A}\cap\mathcal{B}$ will be called the \textit{coheart} of the co-$t$-structure.}
\edf

It would be hoped that the coheart of a co-$t$-structure on $\T$ would be an abelian subcategory of $\T$. Unfortunately, this is not the case, see \cite{Bondarko}. However, in this paper we present an example in which the coheart does turn out to be an abelian category.

In the next section we give some examples of co-$t$-structures.

\section{Some examples of $t$-structures and co-$t$-structures}\label{examples}

\subsection{A canonical example} 

The following example is taken from \cite{Bondarko}. Let $\mathcal{C}$ be an additive category and let $\mathcal{K}(\mathcal{C})$ be its homotopy category. We claim that the following pair of full subcategories of $\mathcal{K}(\mathcal{C})$ forms a co-$t$-structure on $\T=\mathcal{K}(\mathcal{C})$. Let
\begin{eqnarray*}
\A & = & \{\textnormal{complexes in $\mathcal{K}(\mathcal{C})$ isomorphic to complexes } C \ | \ C^{i}=0 \textnormal{ for } i<0\}, \\
\mathcal{B} & = & \{\textnormal{complexes in $\mathcal{K}(\mathcal{C})$ isomorphic to complexes } C \ | \ C^{i}=0 \textnormal{ for } i>0\}.
\end{eqnarray*}
It is clear that $\A$ and $\mathcal{B}$ are closed under direct summands, and that $\Sigma^{-1}\A\subseteq \A$ and $\mathcal{B}\subseteq\Sigma^{-1}\mathcal{B}$. It is also clear that $\Hom{}{\Sigma^{-1}\A}{\mathcal{B}}=0$. We need to show that property $(3)$ of Definition \ref{cot} holds. Suppose $X$ is an object of $\mathcal{K}(\mathcal{C})$:
$$\xymatrix{X: & \cdots\ar[r] & X^{-2}\ar[r] & X^{-1}\ar[r] & X^{0}\ar[r] & X^{1}\ar[r] & X^{2}\ar[r] & \cdots.}$$
We obtain the following semi-split short exact sequence of complexes:
$$\xymatrix{\Sigma^{-1}A:\ar[d] & \cdots\ar[r] & 0\ar[r]\ar[d] & 0\ar[r]\ar[d] & 0\ar[r]\ar[d] & X^{1}\ar[r]\ar@{=}[d] & X^{2}\ar[r]\ar@{=}[d] & \cdots \\
X:\ar[d] & \cdots\ar[r] & X^{-2}\ar[r]\ar@{=}[d] & X^{-1}\ar[r]\ar@{=}[d] & X^{0}\ar[r]\ar@{=}[d] & X^{1}\ar[r]\ar[d] & X^{2}\ar[r]\ar[d] & \cdots \\
B: & \cdots\ar[r] & X^{-2}\ar[r] & X^{-1}\ar[r] & X^{0}\ar[r] & 0\ar[r] & 0\ar[r] & \cdots }$$
which gives us a distinguished triangle
$$\backtri{A}{X}{B}$$
in $\mathcal{K}(\mathcal{C})$. Hence $(\A,\mathcal{B})$ is a co-$t$-structure on $\T=\mathcal{K}(\mathcal{C})$. Moreover, it is non-degenerate and its coheart is just the class of complexes sitting in degree zero.

\subsection{Chain and cochain DGAs}

Recall from the introduction that a DGA $R$ is called a \textit{chain DGA} if $H^{i}(R)=0$ for all $i>0$; similarly, a DGA $R$ is called a \textit{cochain DGA} if $H^{i}(R)=0$ for all $i<0$. A cochain DGA $R$ is called \textit{simply connected} if, in addition, $H^{0}(R)$ is a division ring and $H^{1}(R)=0$.

\eg
{Let $R$ be a chain DGA. Let $\D(R)$ be the derived category of DG $R$-modules, see \cite{Bernstein}, and define a pair of subcategories of $\D(R)$ as follows:
\begin{eqnarray*}
\mathcal{X} & = & \{M\in\D(R) \ | \ H^{i}(M)=0 \textnormal{ for } i>0\}, \\
\mathcal{Y} & = & \{M\in\D(R) \ | \ H^{i}(M)=0 \textnormal{ for } i<0\}.
\end{eqnarray*}
It is easy to show that the pair $(\mathcal{X},\mathcal{Y})$ forms a $t$-structure on $\D(R)$. Again, this $t$-structure is non-degenerate and its heart consists of DG $R$-modules whose cohomology is concentrated in degree zero.}
\label{chain1}
\eeg

\eg
{Let $R$ be a simply connected cochain DGA. Let $\D(R)$ be the derived category of DG $R$-modules and define a pair of subcategories of $\D(R)$ as follows:
\begin{eqnarray*}
\A & = & \{M\in\D(R) \ | \ H^{i}(M)=0 \textnormal{ for } i<0\}, \\
\mathcal{B} & = & \{M\in\D(R) \ | \ H^{i}(M)=0 \textnormal{ for } i>0\}.
\end{eqnarray*}
It is easy to show that the pair $(\A,\mathcal{B})$ forms a co-$t$-structure on $\D(R)$. As in Example \ref{chain1}, this co-$t$-structure is non-degenerate and its coheart consists of DG $R$-modules whose cohomology sits in degree zero.}
\label{cochain1}
\eeg

Simply connected cochain DGAs arise naturally in algebraic topology as the cochain algebras of simply connected CW-complexes, see \cite{Felix} and \cite{Rotman}, for example.

\subsection{A $t$-structure obtained from a rigid object}

Example \ref{chain1} can be abstracted to an arbitrary triangulated category by looking at objects behaving like chain DGAs. Let $R$ be a chain DGA, and recall from the introduction that, given a DG $R$-module $M$ we have
$$H^{i}(M)=\Hom{\D(R)}{R}{\Sigma^{i}M} \textnormal{ for } i\in\mathbb{Z}.$$
Thus, $R$ being a chain DGA means that $\Hom{\D(R)}{R}{\Sigma^{i}R}=0$ for $i>0$.

Now let $\T$ be an arbitrary triangulated category with set indexed coproducts. In the introduction, an object $S$ of $\T$ was called \textit{rigid} if we had
$$\Hom{\T}{S}{\Sigma^{i}S}=0 \textnormal{ for } i>0.$$
Hence, a DGA $R$ is a chain DGA if an only if it is a rigid object in its derived category $\D(R)$. If we replace the chain DGA $R$ with some suitably nice rigid object $S$ of $\T$, the following is a candidate for a $t$-structure on $\T$:
\begin{eqnarray*}
\mathcal{X} & = & \{X\in\T \ | \ \Hom{\T}{S}{\Sigma^{i}X}=0 \textnormal{ for } i>0\}, \\
\mathcal{Y} & = & \{X\in\T\ | \ \Hom{\T}{S}{\Sigma^{i}X}=0 \textnormal{ for } i<0\}.
\end{eqnarray*}
The suitably nice conditions we must place on $S$ to obtain this $t$-structure are that $S$ must be a compact object of $\T$ and the set $\{\Sigma^{i}S \ | \ i\in\mathbb{Z}\}$ must be a generating set for $\T$. Before stating the theorem in full, we recall the notions of a compact object and a generating set.

An object $S$ in a triangulated category $\T$ with set indexed coproducts is \textit{compact} if the functor $\Hom{}{S}{-}$ commutes with set indexed coproducts, that is the canonical map is an isomorphism
$$\Hom{}{S}{\coprod_{i\in I}X_{i}}\cong \coprod_{i\in I}\Hom{}{S}{X_{i}}$$
for all families of objects $\{X_{i}\}_{i\in I}$ of $\T$ indexed by a set $I$; see \cite{Neeman} and \cite{Neeman2}. A DGA $R$ is trivially a compact object of $\D(R)$.

A set of objects $\mathcal{G}$ in a triangulated category $\T$ is called a \textit{generating set} for $\T$ if given any object $X$ of $\T$ with $\Hom{}{G}{X}=0$ for all objects $G$ of $\mathcal{G}$, we have $X=0$.

Example \ref{chain1} is a special case of the following theorem of Hoshino, Kato and Miyachi, which appears in \cite{Miyachi}.

\thm
[\cite{Miyachi}, Theorem 1.3]
Let $\T$ be a triangulated category with set indexed coproducts. Suppose $S$ is a compact rigid object of $\T$ and assume that $\{\Sigma^{i}S \ | \ i\in\mathbb{Z}\}$ is a generating set for $\T$. Then the following forms a non-degenerate $t$-structure on $\T$:
\begin{eqnarray*}
\mathcal{X} & = & \{X\in\T \ | \ \Hom{\T}{S}{\Sigma^{i}X}=0 \textnormal{ for } i>0\}, \\
\mathcal{Y} & = & \{X\in\T\ | \ \Hom{\T}{S}{\Sigma^{i}X}=0 \textnormal{ for } i<0\}.
\end{eqnarray*}
Moreover, its heart $\mathcal{H}=\mathcal{X}\cap\mathcal{Y}$ is an admissible abelian subcategory of $\T$ in the sense of \cite{Beilinson}, and the functor
$$\Hom{}{S}{-}:\mathcal{H}\rightarrow\textnormal{Mod}(\textnormal{End}(S)^{\textnormal{op}})$$
is an equivalence of categories.
\label{miyachitheorem}
\ethm

The $t$-structure is induced as follows: suppose $X$ is an object of $\T$, a morphism $\alpha:X\rightarrow Y$ with $Y\in \mathcal{Y}$ is constructed such that, given any other object $Y'\in \mathcal{Y}$ and a morphism  $X\rightarrow Y'$, then this morphism factors uniquely through $\alpha:X\rightarrow Y$,
$$\xymatrix{X\ar[r]^-{\alpha}\ar[dr] & Y\ar@{-->}[d]^-{\exists !} \\
& Y'.}$$

\subsection{A co-$t$-structure obtained from a corigid object}

In this paper we shall look at the structure which is induced by an object behaving like a cochain DGA. The subsequent sections of this paper are devoted to proving the following theorem, which is the cochain analogue, or dual, of Theorem \ref{miyachitheorem}.

\thm
Let $\T$ be a triangulated category with set indexed coproducts. Suppose $S$ is a compact simply connected corigid object of $\T$ and assume that $\{\Sigma^{i}S \ | \ i\in\mathbb{Z}\}$ is a generating set for $\T$. Then the following forms a non-degenerate co-$t$-structure on $\T$:
\begin{eqnarray*}
\mathcal{A} & = & \{X\in\T \ | \ \Hom{\T}{S}{\Sigma^{i}X}=0 \textnormal{ for } i<0\}, \\
\mathcal{B} & = & \{X\in\T\ | \ \Hom{\T}{S}{\Sigma^{i}X}=0 \textnormal{ for } i>0\}.
\end{eqnarray*}
Moreover, its coheart $\mathcal{C}=\mathcal{A}\cap\mathcal{B}$ is an abelian subcategory of $\T$, and the functor
$$\Hom{}{S}{-}:\mathcal{C}\rightarrow\textnormal{Mod}(\textnormal{End}(S)^{\textnormal{op}})$$
is an equivalence of categories.
\label{cochaintheorem}
\ethm

Theorem \ref{cochaintheorem} is the natural generalisation of Example \ref{cochain1} in the same way that Theorem \ref{miyachitheorem} is the natural generalisation of Example \ref{chain1}.

In Section \ref{preenvelopes} we shall show that given any object $M$ of $\T$, there exists a morphism $\mu:M\rightarrow\bar{M}$ with $\bar{M}\in S^{\perp_{\infty}}$ such that, given any other object $N\in S^{\perp_{\infty}}$ and a morphism $M\rightarrow N$, then this morphism factors through $\mu:M\rightarrow\bar{M}$,
$$\xymatrix{M\ar[r]^-{\alpha}\ar[dr] & \bar{M}\ar@{-->}[d]^-{\exists} \\
& N.}$$ 
However, the factorisation is not necessarily unique. Thus we obtain an $S^{\perp_{\infty}}$-preenvelope. Note that in Theorem \ref{cochaintheorem} above, $\mathcal{B}=S^{\perp_{\infty}}$. In Section \ref{particular} we show that this $S^{\perp_{\infty}}$-preenvelope induces a non-degenerate co-$t$-structure on $\T$, while Section \ref{coheart} is dedicated to proving that the coheart of this non-degenerate co-$t$-structure is equivalent to the module category $\textnormal{Mod}(\textnormal{End}(S)^{\textnormal{op}})$, and hence abelian.

We now make precise what we mean by an object of a triangulated category behaving like a cochain DGA. Following the definition of an $n$-rigid object in a triangulated category of Iyama and Yoshino in \cite{Iyama}, we make the following definitions of an $n$-corigid object and a corigid object.

\df
{An object $S$ of $\T$ will be called \textit{$n$-corigid} if we have $$\Hom{}{\Sigma^{i}S}{S}=0 \textnormal{ for } 0<i<n.$$ An object $S$ of $\T$ will be called \textit{corigid} if we have $$\Hom{}{\Sigma^{i}S}{S}=0 \textnormal{ for } i>0.$$}
\edf

Note that a DGA $R$ is a cochain DGA if and only if it is a corigid object in its derived category $\D(R)$.

\df
{Let $S$ be an object of $\T$. We shall call $S$ a \textit{simply connected corigid object} of $\T$ if it satisfies the following assumptions:}

\textnormal{$(1)$ $S$ is corigid, that is, $\Hom{}{\Sigma^{i}S}{S}= 0$ for $i>0$;}

\textnormal{$(2)$ $\Hom{}{S}{\Sigma S}=0$;}

\textnormal{$(3)$ $\Endo{S}$ is a division ring.}
\label{cochain}
\edf

Note that a DGA $R$ is a  simply connected cochain DGA if and only if it is a simply connected corigid object in its derived category $\D(R)$.

\rem
{For technical reasons, in the cochain analogue of Theorem \ref{miyachitheorem} we must also insist that $S$ is simply connected in the sense of Definition \ref{cochain}. This is due to the lack of symmetry in the theory of chain and cochain DGAs mentioned in the introduction: one is able to construct a theory for chain DGAs, but in order to construct a viable dual theory for cochain DGAs one has to introduce the assumption of simply connectedness; see, for example, \cite{Avramov}.}
\label{dichotomy}
\erem

\section{Existence of an $S^{\perp_{\infty}}$-preenvelope}\label{preenvelopes}

In order to obtain an $S^{\perp_{\infty}}$-preenvelope, we first show how to construct an $S^{\perp_{n}}$-preenvelope for each $n\in\mathbb{N}$. It is useful to refer to a simply connected $n$-corigid object of a triangulated category:

\df
{Let $S$ be an object of $\T$. We shall call $S$ a \textit{simply connected $n$-corigid object} of $\T$ if it satisfies the following assumptions:}

\textnormal{$(1)$ $S$ is $n$-corigid, that is, $\Hom{}{\Sigma^{i}S}{S}= 0$ for $0<i<n$;}

\textnormal{$(2)$ $\Hom{}{S}{\Sigma S}=0$;}

\textnormal{$(3)$ $\Endo{S}$ is a division ring.}
\label{n_cochain}
\edf

\prop
Let $\T$ be a triangulated category with set indexed coproducts. Suppose $S$ is a compact simply connected $(n+1)$-corigid object of $\T$. Then, for each object $M$ of $\T$ there exists an  $S^{\perp_{n}}$-preenvelope $\mu :M\rightarrow \bar{M}$.
\label{prop2}
\eprop

\pf
Let $M$ be an arbitrary object of $\T$. We first construct a chain of objects and morphisms,
$$M=M_{0}\rightlabel{\mu_{0}} M_{1}\rightlabel{\mu_{1}} M_{2}\rightlabel{\mu_{2}} M_{3}\rightlabel{\mu_{3}} \cdots \rightlabel{\mu_{n-1}} M_{n}=\bar{M}$$
with $M_{k}\in S^{\perp_{k}}$ for each $k\geqslant 1$, inductively using distinguished triangles. Secondly, we verify that the composite of these maps is an $S^{\perp_{n}}$-preenvelope.

Write $M=M_{0}$. Let $n=1$; we  construct an object $M_{1}$ and a morphism $\mu_{0}:M_{0}\rightarrow M_{1}$ such that $\Hom{}{S}{\Sigma M_{1}}=0$. If $\Hom{}{S}{\Sigma M_{0}}=0$ then set $M_{1}=M_{0}$ and $\mu_{0}=1_{M_{0}}$, the identity map on $M_{0}$. If not, we can choose a, possibly infinite, coproduct $S^{(m_{1})}$ of copies of $S$ and a nonzero morphism $S^{(m_{1})}\rightarrow \Sigma M_{0}$ which becomes a surjection under the functor $\Hom{}{S}{-}$. Since the endomorphism ring $\Endo{S}$ is a division ring we can, moreover, choose $m_{1}$ so that this morphism becomes an isomorphism under $\Hom{}{S}{-}$. We now extend this morphism to a distinguished triangle:
\begin{equation}
\tri{S^{(m_{1})}}{\Sigma M_{0}}{\Sigma M_{1}}.
\label{star}
\end{equation}
Applying $\Hom{}{S}{-}$ to \eqref{star} gives the exact sequence:
$$\Hom{}{S}{S^{(m_{1})}}\rightiso \Hom{}{S}{\Sigma M_{0}} \rightarrow \Hom{}{S}{\Sigma M_{1}}\rightarrow \Hom{}{S}{\Sigma S^{(m_{1})}}.$$
Since $\Hom{}{S}{\Sigma S^{(m_{1})}}=0$, we get $\Hom{}{S}{\Sigma M_{1}}=0$.

Now suppose $k\geqslant 1$ and suppose we have constructed a chain of objects and morphisms
$$M=M_{0}\rightlabel{\mu_{0}} M_{1}\rightlabel{\mu_{1}} M_{2}\rightlabel{\mu_{2}} M_{3}\rightlabel{\mu_{3}} \cdots \rightlabel{\mu_{k-1}} M_{k}$$
with $M_{i}\in S^{\perp_{i}}$ for $1\leqslant i\leqslant k$, and where $\mu_{i}:M_{i}\rightarrow M_{i+1}$ is either the identity map or sits in a distinguished triangle
\begin{equation*}
\Sigma^{-(i+1)}S^{(m_{i+1})}\rightarrow M_{i}\rightlabel{\mu_{i}} M_{i+1}\rightarrow \Sigma^{-i}S^{(m_{i+1})}.
\label{triangle_i}
\end{equation*}
If $\Hom{}{S}{\Sigma^{k+1}M_{k+1}}=0$ then set $M_{k+1}=M_{k}$ and take $\mu_{k}:M_{k}\rightarrow M_{k+1}$ to be the identity map $1_{M_{k}}$. If not, we can choose a, possibly infinite, coproduct $S^{(m_{k+1})}$ of copies of $S$ and a nonzero morphism $S^{(m_{k+1})}\rightarrow\Sigma^{k+1}M_{k}$ which becomes an isomorphism under $\Hom{}{S}{-}$, and then extend it to a distinguished triangle:
\begin{equation}
\tri{S^{(m_{k+1})}}{\Sigma^{k+1}M_{k}}{\Sigma^{k+1}M_{k+1}}.
\label{doublestar}
\end{equation}
As above, an argument from the long exact sequence of Hom-sets arising from \eqref{doublestar} shows that $$\Hom{}{S}{\Sigma^{i}M_{k+1}}=0 \textnormal{ for }i=1,\ldots,k+1.$$ The case $i=k$ follows by the injectivity of $\Hom{}{S}{S^{(m_{k+1})}}\rightiso \Hom{}{S}{\Sigma^{k+1}M_{k}}$; and the case $i=k+1$ by its surjectivity.

Hence, inductively we obtain  a chain of objects and morphisms of $\T$,
\begin{equation}
M=M_{0}\rightlabel{\mu_{0}} M_{1}\rightlabel{\mu_{1}} M_{2}\rightlabel{\mu_{2}} M_{3}\rightlabel{\mu_{3}}\cdots\rightlabel{\mu_{n-1}} M_{n},
\label{composite}
\end{equation}
where each map $\mu_{k}:M_{k}\rightarrow M_{k+1}$ is either the identity map or sits in a distinguished triangle
\begin{equation}
\Sigma^{-(k+1)}S^{(m_{k+1})}\rightarrow M_{k}\rightlabel{\mu_{k}} M_{k+1}\rightarrow \Sigma^{-k}S^{(m_{k+1})}.
\label{triangle_k}
\end{equation}
To see that the composite $\mu=\mu_{n-1}\circ\cdots\circ\mu_{1}\circ\mu_{0}$ from \eqref{composite} is an $S^{\perp_{n}}$-preenvelope, we shall show that for each $X\in S^{\perp_{n}}$ the map $\Hom{}{M_{k+1}}{X}\rightarrow \Hom{}{M_{k}}{X}$ induced by $\mu_{k}$ is a surjection.
Without loss of generality we may assume that each map $\mu_{k}$ sits in a distinguished triangle \eqref{triangle_k} above, because if $\mu_{k}=1_{M_{k}}$, then the map $\Hom{}{M_{k+1}}{X}\rightarrow \Hom{}{M_{k}}{X}$ is trivially an isomorphism for all $X\in \T$.

Let $X\in S^{\perp_{n}}$;  applying $\Hom{}{-}{X}$ to distinguished triangle \eqref{triangle_k}, we get the long exact sequence of Hom-sets below:
$$(\Sigma^{-k}S^{(m_{k+1})},X)\rightarrow (M_{k+1},X)\rightarrow (M_{k},X)\rightarrow (\Sigma^{-(k+1)}S^{(m_{k+1})},X),$$
where we have written $(A,B)$ as a  shorthand for $\Hom{}{A}{B}$.
Now since we have  $\Hom{}{S^{(m_{k+1})}}{\Sigma^{k}X} = \Hom{}{S^{(m_{k+1})}}{\Sigma^{(k+1)}X}=0$ for $k=1,\ldots,n-1$, the map
$$\Hom{}{M_{k+1}}{X}\rightarrow \Hom{}{M_{k}}{X}$$
induced by $\mu_{k}$ is an isomorphism for $k=1,\ldots, n-1$ and a surjection for $k=0$. Hence, writing $\bar{M}=M_{n}$, the composite $\mu:M\rightarrow \bar{M}$ is an $S^{\perp_{n}}$-preenvelope.
\epf

\lemma
\label{corollary0}
Suppose further that $S$ is a simply connected corigid object of $\T$. Then, for the $S^{\perp_{n}}$-preenvelope, $\mu:M\rightarrow \bar{M}$, obtained in Proposition \ref{prop2}, we have that $$\Hom{}{S}{\Sigma^{i}\mu}:\Hom{}{S}{\Sigma^{i}M}\rightarrow \Hom{}{S}{\Sigma^{i}\bar{M}}$$ is an isomorphism for all $i<1$.
\elemma

\pf
Applying the functor $\Hom{}{S}{-}$ to distinguished triangle \eqref{triangle_k}, 
$$\Sigma^{-(k+1)}S^{(m_{k+1})}\rightarrow M_{k}\rightlabel{\mu_{k}} M_{k+1}\rightarrow \Sigma^{-k}S^{(m_{k+1})},$$ 
for $0\leqslant k< n$ in the proof of Proposition \ref{prop2} shows that $$\Hom{}{S}{\Sigma^{i}\mu}:\Hom{}{S}{\Sigma^{i}M_{k}}\rightarrow \Hom{}{S}{\Sigma^{i}M_{k+1}}$$ is an isomorphism for all $i<k+1$. The isomorphism for $i=k$ follows by the fact that the morphism $S^{(m_{k+1})}\rightarrow \Sigma^{k+1}M_{k}$ in \eqref{triangle_k} is constructed to be an isomorphism under $\Hom{}{S}{-}$. Hence the composite $\mu=\mu_{n-1}\circ\cdots\circ\mu_{1}\circ\mu_{0}$ is an isomorphism under $\Hom{}{S}{\Sigma^{i}-}$ for all $i<1$.
\epf

\

In order to obtain an $S^{\perp_{\infty}}$-preenvelope we need to introduce the key tool, which is called the \textit{homotopy colimit}. The following definition is taken from \cite{Neeman}.

\df
{Let $\T$ be a triangulated category with set indexed coproducts. Let 
$$X_{0}\rightlabel{f_{0}} X_{1}\rightlabel{f_{1}} X_{2} \rightlabel{f_{2}} X_{3}\rightlabel{f_{3}}\cdots$$ 
be a sequence of objects and morphisms in $\T$. The \textit{homotopy colimit} $\hocolim{X_{i}}$ is constructed by extending the map
$$\coprod_{i=0}^{\infty}X_{i}\rightlabel{1-\rm{shift}}\coprod_{i=0}^{\infty}X_{i}$$
to a distinguished triangle:
$$\trilabel{\coprod_{i=0}^{\infty}X_{i}}{\coprod_{i=0}^{\infty}X_{i}}{\hocolim{X_{i}}}{1-\rm{shift}}{}{}.$$}
\edf

We will need the following lemma.

\lemma
[\cite{Neeman}, Lemma 2.8]
\label{lemma1}
Suppose $S$ is a compact object of a triangulated category $\T$ and we have a sequence of objects and morphisms of $\T$:
$$X_{0}\rightarrow X_{1}\rightarrow X_{2}\rightarrow X_{3}\rightarrow \cdots$$
then $\textnormal{colim}\,(\Hom{}{S}{X_{n}})\cong \Hom{}{S}{\textnormal{hocolim}X_{n}}$.
\elemma

\prop
\label{prop3}
Let $\T$ be a triangulated category with set indexed coproducts. Suppose $S$ is a compact simply connected corigid object of $\T$. Then, for each object $M$ of $\T$ there exists an $S^{\perp_{\infty}}$-preenvelope $\mu:M\rightarrow \bar{M}$.
\eprop

\pf
Let $M$ be an object of $\T$ and write $M=M_{0}$. Let $X\in S^{\perp_{\infty}}$ and suppose we have a morphism $\alpha_{0}:M_{0}\rightarrow X$. By the argument of Proposition \ref{prop2} we can construct the following commutative diagram:
$$\xymatrix{ & X & & & &  \\
M_{0}\ar[r]_-{\mu_{0}}\ar[ur]^-{\alpha_{0}} & M_{1}\ar[r]_-{\mu_{1}}\ar[u]^-{\alpha_{1}} & M_{2}\ar[r]_-{\mu_{2}}\ar[ul]^-{\alpha_{2}} & \cdots\ar[r] & M_{n}\ar[r]_-{\mu_{n}}\ar[ulll]_-{\alpha_{n}} & \cdots \ .}$$
with $M_{n}\in S^{\perp_{n}}$ for each $n\geqslant 1$. We now construct the homotopy colimit, $\hocolim{M_{i}}$.
By construction, the composite $$\coprod_{i=0}^{\infty} M_{i}\rightlabel{1-\textnormal{shift}} \coprod_{i=0}^{\infty} M_{i}\rightlabel{\langle\alpha_{i}\rangle} X$$ is zero, so that we have the following commutative diagram:
$$\xymatrix{\coprod_{i=0}^{\infty} M_{i}\ar[r]^-{1-\textnormal{shift}}\ar[dr]_-{0} & \coprod_{i=0}^{\infty} M_{i}\ar[r]\ar[d]^-{\langle\alpha_{i}\rangle} & \hocolim{M_{i}}\ar[r]\ar@{-->}[dl]^-{\exists} & \Sigma \coprod M_{i} \\
 & X & & }.$$
That is, every morphism $M\rightarrow X$ factors through $\hocolim{M_{i}}\rightarrow X$. 

Now we have:
\begin{eqnarray*}
\Hom{}{S}{\Sigma^{j}\hocolim{M_{i}}} & \cong & \Hom{}{S}{\hocolim{\Sigma^{j}M_{i}}} \\
				     & \cong & \textnormal{colim} \ \Hom{}{S}{\Sigma^{j}M_{i}} \\
				     & =     & 0
\end{eqnarray*}
for $j\geqslant 1$. We obtain the first isomorphism because the homotopy colimit commutes with the suspension functor and the second isomorphism by Lemma \ref{lemma1}. The final equality is a consequence of the fact that $\Hom{}{S}{\Sigma^{j}M_{i}}=0$ for $i$ sufficiently large and $j\geqslant 1$. Hence we have $\hocolim{M_{i}}\in S^{\perp_{\infty}}$. Therefore, setting $\bar{M}=\hocolim{M_{i}}$, we obtain an $S^{\perp_{\infty}}$-preenvelope $\mu:M\rightarrow \bar{M}$.
\epf

\section{A co-$t$-structure induced by a compact simply connected corigid object}\label{particular}

The aim of this section is to give a proof of the following theorem, which is the first half of Theorem \ref{cochaintheorem}.

\thm
\label{prop4}
Let $\T$ be a triangulated category with set indexed coproducts. Suppose $S$ is a compact simply connected corigid  object of $\T$. Further assume that $\{\Sigma^{i}S \ | \ i\in\Z\}$ is a generating set in $\T$. Then the following forms a non-degenerate co-$t$-structure on $\T$:
\begin{eqnarray*}
\A & = & \{X\in \T \ | \ \Hom{}{S}{\Sigma^{i}X}=0 \textnormal{ for } i<0\}, \\
\mathcal{B} & = & \{X\in \T \ | \ \Hom{}{S}{\Sigma^{i}X}=0 \textnormal{ for } i>0\}. 
\end{eqnarray*}
Note that $\mathcal{B}=S^{\perp_{\infty}}$.
\ethm

In order to prove this we need a lemma analogous to Lemma \ref{corollary0}. This is an immediate consequence of the next lemma.

\lemma
\label{third}
Let $S$ be a compact object of a triangulated category $\T$ and suppose we have a sequence of objects and morphisms $$X_{0}\rightlabel{\alpha_{0}} X_{1}\rightlabel{\alpha_{1}} X_{2}\rightlabel{\alpha_{2}} X_{3}\rightlabel{\alpha_{3}} \cdots$$
such that $\Hom{}{S}{\alpha_{n}}:\Hom{}{S}{X_{n}}\rightarrow \Hom{}{S}{X_{n+1}}$ is an isomorphism for each $n\geqslant 0$. Then $\Hom{}{S}{X_{0}}\cong \Hom{}{S}{\hocolim{X_{n}}}$.
\elemma

\pf
It is well-known that the filtered colimit, $\textnormal{colim}\, \Hom{}{S}{X_{n}}$, is isomorphic to $\Hom{}{S}{X_{0}}$. The assertion now follows by Lemma \ref{lemma1}.
\epf

\lemma
\label{corollary1}
Under the assumptions of Proposition \ref{prop3} we have that  $$\Hom{}{S}{\Sigma^{i}\mu}:\Hom{}{S}{\Sigma^{i}M}\rightarrow\Hom{}{S}{\Sigma^{i}\bar{M}}$$ is an isomorphism for $i<1$.
\elemma

\pf
In Lemma \ref{corollary0}, $\Hom{}{S}{\Sigma^{i}\mu}:\Hom{}{S}{\Sigma^{i}M_{k}}\rightarrow \Hom{}{S}{\Sigma^{i}M_{k+1}}$ is an isomorphism for $i<1$. Now apply Lemma \ref{third}.
\epf

\

\noindent\textbf{Proof of Theorem \ref{prop4}:} Conditions $(0)$ and $(1)$ of the definition of a co-$t$-structure are clear.

In order to show $(2)$ assume $X\in\Sigma^{-1}A$ and $Y\in\mathcal{B}$. Recall that $\mathcal{B}=S^{\perp_{\infty}}$. By Proposition \ref{prop3} there exists  an $S^{\perp_{\infty}}$-preenvelope $\mu:X\rightarrow\bar{X}$, that is, we have a surjection of Hom spaces
$$\Hom{}{\bar{X}}{Y}\twoheadrightarrow \Hom{}{X}{Y}.$$ It is therefore sufficient to show $\Hom{}{\bar{X}}{Y}=0$. 

We have the following isomorphism of Hom-spaces and trivial Hom-spaces:
\begin{eqnarray*}
\Hom{}{S}{\Sigma^{i}X} & \cong & \Hom{}{S}{\Sigma^{i}\bar{X}} \textnormal{ for all } i<1 \textnormal{\ \ \ (Lemma \ref{corollary1})} \\
\Hom{}{S}{\Sigma^{i}X} & = & 0 \textnormal{ for all } i<1 \textnormal{\ \ \ (since $X\in\Sigma^{-1}\mathcal{A}$)} \\
\Hom{}{S}{\Sigma^{i}\bar{X}} & = & 0 \textnormal{ for all } i>0 \textnormal{\ \ \ (since $\bar{X}\in\mathcal{B}$)}.
\end{eqnarray*}
It follows that $\Hom{}{S}{\Sigma^{i}\bar{X}}=0$ for all $i\in\mathbb{Z}$. The assumption that $\{\Sigma^{i}S \ | \  i\in\mathbb{Z}\}$ is a generating set in $\T$ implies that $\bar{X}=0$. Thus $\Hom{}{\bar{X}}{Y}=0$ and we see that $\Hom{}{X}{Y}=0$. Hence $\Hom{}{\Sigma^{-1}\mathcal{A}}{\mathcal{B}}=0$.

We next show condition $(3)$. Suppose $X$ is an object of $\T$. By Proposition \ref{prop3} there is an $S^{\perp_{\infty}}$-preenvelope $\mu:X\rightarrow \bar{X}$. Write $B=\bar{X}$ and extend the morphism $\mu:X\rightarrow B$ to a distinguished triangle:
\begin{eqnarray}
\backtri{A}{X}{B}.
\label{triangle1}
\end{eqnarray}
We claim that $\Hom{}{S}{\Sigma^{i}A}=0$ for $i<0$. Consider the following long exact sequence obtained from (\ref{triangle1}):
$$\Hom{}{S}{\Sigma^{i-1}A}\rightarrow \Hom{}{S}{\Sigma^{i}X}\rightarrow \Hom{}{S}{\Sigma^{i}B}\rightarrow \Hom{}{S}{\Sigma^{i}A}.$$
Now, by Lemma \ref{corollary1}, we see that $\Hom{}{S}{\Sigma^{i}X}\rightarrow\Hom{}{S}{\Sigma^{i}B}$ is an isomorphism for all $i<1$. Hence $\Hom{}{S}{\Sigma^{i}A}=0$ for all $i<0$ and $A\in \mathcal{A}$. Hence the distinguished triangle in (\ref{triangle1}) above gives us the required distinguished triangle.

It is clear that $\cap_{n\in\Z}\Sigma^{n}\A=\cap_{n\in\Z}\Sigma^{n}\mathcal{B}=\{0\}$ because $\{\Sigma^{i}S \ | \ i\in\Z\}$ is a generating set for $\T$. Hence $(\A,\mathcal{B})$ is a non-degenerate co-$t$-structure on $\T$. $\Box$

\rem
{In \cite{Bondarko}, a co-$t$-structure $(\A,\mathcal{B})$ is called \textit{right adjacent} to a $t$-structure $(\mathcal{X},\mathcal{Y})$ if $\A=\mathcal{Y}$. By \cite{Beligiannis_Reiten}, the full subcategory $\A$ of $\T$ occurs in a $t$-structure $(\mathcal{X},\A)$ on $\T$, where $$\mathcal{X}={}^{\perp}\A : = \{X\in\T \, |\, \Hom{}{X}{A}=0 \textnormal{ for all } A\in\A\}.$$ Therefore, the co-$t$-structure on $\T$ obtained in Theorem \ref{prop4} is right adjacent to the $t$-structure $(\mathcal{X},\A)$.}
\erem

\section{The coheart of the co-$t$-structure of Theorem \ref{prop4}}\label{coheart}

In Theorem \ref{miyachitheorem}, Hoshino, Kato and Miyachi not only obtain a non-degenerate $t$-structure on a triangulated category $\T$, but they also show that its heart, which is admissible abelian, is equivalent to the module category $\textnormal{Mod}(\textnormal{End}(S)^{\textnormal{op}})$.  We shall show that the coheart of the co-$t$-structure obtained in Theorem \ref{prop4} is equivalent to $\textnormal{Mod}(\Endo{S}^{\textnormal{op}})$, where $S$ is the object from Theorem \ref{prop4} and where $\Endo{S}^{\textnormal{op}}$ indicates that this is the category of right $\Endo{S}$-modules. This is the second half of Theorem \ref{cochaintheorem}, and will complete the proof of the cochain analogue of Theorem \ref{miyachitheorem}.

\setup
{Throughout this section, we shall consider the co-$t$-structure of Theorem \ref{prop4}; that is, let $\T$ be a triangulated category with set indexed coproducts, suppose $S$ is a compact simply connected corigid object of $\T$. Furthermore,  assume that $\{\Sigma^{i}S \ | \ i\in\mathbb{Z}\}$ is a generating set for $\T$. Then, by Theorem \ref{prop4}, the following is a non-degenerate co-$t$-structure on $\T$:
\begin{eqnarray*}
\A & = & \{X\in \T \ | \ \Hom{}{S}{\Sigma^{i}X}=0 \textnormal{ for } i<0\}, \\
\mathcal{B} & = & \{X\in \T \ | \ \Hom{}{S}{\Sigma^{i}X}=0 \textnormal{ for } i>0\}. 
\end{eqnarray*}
Let $\mathcal{C}=\mathcal{A}\cap\mathcal{B}$ be the coheart of this co-$t$-structure.}
\label{setup}
\esetup

\lemma
\label{surjects}
Under the hypotheses of Setup \ref{setup}  the functor $$\Hom{}{S}{-}:\mathcal{C}\rightarrow \textnormal{Mod}(k^\textnormal{op}),$$ where $k=\Endo{S}$, is dense.
\elemma

\pf
Consider the object $S$ of $\T$. Since $(\mathcal{A},\mathcal{B})$ forms a co-$t$-structure on $\T$, there is a distinguished triangle
\begin{equation}
\backtrilabel{A}{S}{B}{\alpha}{}{}
\label{triangle4}
\end{equation}
with $A\in\mathcal{A}$ and $B\in\mathcal{B}$. Applying the functor $\Hom{}{S}{-}$ to \eqref{triangle4} gives the following long exact sequence:
\begin{eqnarray}
\Hom{}{S}{\Sigma^{i}S}\rightarrow \Hom{}{S}{\Sigma^{i}B}\rightarrow \Hom{}{S}{\Sigma^{i}A}.
\label{longexact1}
\end{eqnarray}
In (\ref{longexact1}) we have $\Hom{}{S}{\Sigma^{i}S}=\Hom{}{S}{\Sigma^{i}A}=0$ for all $i<0$ since $A\in\mathcal{A}$, so that $\Hom{}{S}{\Sigma^{i}B}=0$ for all $i<0$. We know that $\Hom{}{S}{\Sigma^{i}B}=0$ for all $i>0$ since $B\in\mathcal{B}$. Therefore, $\Hom{}{S}{\Sigma^{i}B}=0$ for all $i\neq 0$. By Lemma \ref{corollary1}, $\Hom{}{S}{S}\rightarrow \Hom{}{S}{B}$ is an isomorphism. Hence we have:
\begin{displaymath}
\Hom{}{S}{\Sigma^{i}B^{(m)}} = \left\{ \begin{array}{ll}
0 & \textrm{if $i\neq 0$} \\
k^{(m)} & \textrm{if $i=0$.}
\end{array} \right.
\end{displaymath}
where $k=\Endo{S}$. Hence $\Hom{}{S}{-}:\mathcal{C}\rightarrow \textnormal{Mod}(k^{\textnormal{op}})$ is dense.
\epf

\

Note that the fact that $\Hom{}{S}{S}\rightarrow \Hom{}{S}{B}$ is an isomorphism forces $\Hom{}{S}{A}=0$. Therefore, $A$ is an object of $\Sigma^{-1}\mathcal{A}$ rather than just an object of $\mathcal{A}$.

The object $B$ introduced in triangle (\ref{triangle4}) above  has the useful property that any object in the coheart $\mathcal{C}$ can be described in terms of it; we show this in the next lemma.

\lemma
\label{objectdescribe}
Under the hypotheses of Setup \ref{setup} we have that each $M\in\mathcal{C}$ is isomorphic to $B^{(m)}$ for some $m$.
\elemma

\pf
Consider the distinguished triangle \eqref{triangle4} from Lemma \ref{surjects}:
\begin{eqnarray*}
\Sigma^{-1}A\rightarrow S\rightarrow B\rightarrow A
\label{triangle5}
\end{eqnarray*}
with $A\in\Sigma^{-1}\mathcal{A}$ and $B\in\mathcal{C}$. Let $M\in\mathcal{C}$; since $k=\Hom{}{S}{S}$ is a skew field, we can choose $m$ such that $S^{(m)}\rightarrow M$ becomes an isomorphism under $\Hom{}{S}{-}$. Again, by Lemma \ref{corollary1}, the morphism $S\rightarrow B$ becomes an isomorphism under $\Hom{}{S}{-}$.

We may now apply the functor $\Hom{}{-}{M}$ to \eqref{triangle4} and from the long exact sequence notice that the morphism $\Hom{}{B}{M}\rightarrow \Hom{}{S}{M}$ is an isomorphism. So we obtain the commutative diagram:
$$\xymatrix{\Sigma^{-1}A^{(m)}\ar[r]\ar[dr]_-{0} & S^{(m)}\ar[r]\ar[d] & B^{(m)}\ar[r]\ar@{-->}[dl]^-{\exists !} & A^{(m)} \\
 & M & & }$$
where both $S^{(m)}\rightarrow M$ and $S^{(m)}\rightarrow B^{(m)}$ are isomorphisms under $\Hom{}{S}{-}$. Hence the unique map $B^{(m)}\rightarrow M$ making the diagram above commute becomes an isomorphism under $\Hom{}{S}{-}$.

Now extend this unique map $B^{(m)}\rightarrow M$ to a distinguished triangle
\begin{eqnarray}
\tri{B^{(m)}}{M}{Z}
\label{triangle6}
\end{eqnarray}
and apply the functor $\Hom{}{S}{-}$ to give a long exact sequence. One easily sees from this long exact sequence that $\Hom{}{S}{\Sigma^{i}Z}=0$ for $i\neq 0$. The fact that the morphism $B^{(m)}\rightarrow M$ becomes the isomorphism $\Hom{}{S}{B^{(m)}}\rightiso \Hom{}{S}{M}$ forces $\Hom{}{S}{Z}=0$ so that $\Hom{}{S}{\Sigma^{i}Z}=0$ for all $i\in\mathbb{Z}$. Since $\{\Sigma^{i}S \ | \ i\in\mathbb{Z}\}$ is a generating set for $\T$, it follows that $Z=0$. Hence $B^{(m)}\rightarrow M$ is an isomorphism.
\epf

\prop
\label{fullyfaithful}
Under the hypotheses of Setup \ref{setup}, the functor $$\Hom{}{S}{-}:\mathcal{C}\rightarrow \textnormal{Mod}(k^\textnormal{op})$$ is full and faithful.
\eprop

\pf
We first show that $\Hom{}{S}{-}$ is faithful. Again, consider distinguished triangle \eqref{triangle4}: 
\begin{eqnarray*}
\Sigma^{-1}A\rightarrow S\rightarrow B\rightarrow A.
\label{triangle5}
\end{eqnarray*}
By Lemma \ref{objectdescribe}, any object $M\in\mathcal{C}$ is isomorphic to some coproduct $B^{(I)}$, where $I$ is an indexing set and $B^{(I)}$ denotes the, possibly infinite, coproduct $\coprod_{i\in I}B$. Hence, to show fidelity we can consider a morphism $B^{(I)}\rightarrow B^{(J)}$ which becomes zero under $\Hom{}{S}{-}$ and show that it is itself necessarily zero. It is sufficient to show that the composite $B\hookrightarrow B^{(I)}\rightarrow B^{(J)}$ is zero, where the morphism $B\hookrightarrow B^{(I)}$ is just the coproduct inclusion into the $i^{\textnormal{th}}$-summand for each $i\in I$. This puts us in the following situation:
$$\xymatrix{\Sigma^{-1}A\ar[r] & S\ar[r]\ar[ddr]_-{0} & B\ar[r]\ar@{^{(}->}[d] & A\ar@{-->}[ddl]^-{\exists} \\
& & B^{(I)}\ar[d] & \\
& & B^{(J)}. & }$$
But, $A\in\Sigma^{-1}\mathcal{A}$ and $B^{(J)}\in\mathcal{C}=\mathcal{A}\cap\mathcal{B}$, so that $\Hom{}{A}{B^{(J)}}=0$. Therefore, the dotted arrow above is necessarily zero. Hence the composite  $B\hookrightarrow B^{(I)}\rightarrow B^{(J)}$ is zero, showing that $\Hom{}{S}{-}$ is faithful.

We must also show that $\Hom{}{S}{-}$ is full. Suppose we have a morphism
$$\theta: \Hom{}{S}{B^{(I)}}\rightarrow \Hom{}{S}{B^{(J)}},$$
where $I$ and $J$ are again indexing sets. We must construct a morphism
$B^{(I)}\rightarrow B^{(J)}$
which induces $\theta$ under $\Hom{}{S}{-}$. We recall distinguished triangle (\ref{triangle4}) again:
$$\Sigma^{-1}A\longrightarrow S\rightlabel{\sigma} B\longrightarrow A.$$
Note that $\sigma:S\rightarrow B$ becomes an isomorphism under $\Hom{}{-}{B^{(I)}}$ because $B^{(I)}\in\mathcal{C}$. Hence we get the following commutative diagram:
\begin{equation}
\xymatrix{\Hom{}{S}{B^{(I)}}\ar[r]^-{\theta} & \Hom{}{S}{B^{(J)}} \\
\Hom{}{B}{B^{(I)}}\ar[r]_-{\phi}\ar[u]^-{\Hom{}{\sigma}{B^{(I)}}}_-{\sim} & \Hom{}{B}{B^{(J)}}\ar[u]_-{\Hom{}{\sigma}{B^{(J)}}}^-{\sim}}
\label{comm_square}
\end{equation}
where $\phi = \Hom{}{\sigma}{B^{(J)}}^{-1}\circ\theta\circ\Hom{}{\sigma}{B^{(I)}}$. Let $q_{i}:B\hookrightarrow B^{(I)}$ be the $i^{\textrm{th}}$-inclusion of the coproduct and consider its image $\phi(q_{i}):B\rightarrow B^{(J)}$. By the universal property of the coproduct there exists a unique map $\langle\phi(q_{i})\rangle:B^{(I)}\rightarrow B^{(J)}$ such that the following diagram commutes for each $i\in I$:
$$\xymatrix{B\ar[r]^-{q_{i}}\ar[dr]_-{\phi(q_{i})} & B^{(I)}\ar[d]^-{\langle\phi(q_{i})\rangle} \\
& B^{(J)}.}$$

Let us show that $\langle\phi(q_{i})\rangle$ induces $\theta$ under $\Hom{}{S}{-}$. The map $\Hom{}{S}{\sigma}:\Hom{}{S}{S}\rightarrow \Hom{}{S}{B}$ is an isomorphism, therefore, it takes a set of generators for $\Hom{}{S}{S}$ to a set of generators for $\Hom{}{S}{B}$. The vector space $\Hom{}{S}{S}$ is one-dimensional and generated by the identity map on $S$, $1_{S}$, whose image under $\Hom{}{S}{\sigma}$ is $\sigma: S\rightarrow B$. Hence $\Hom{}{S}{B}$ is generated by $\sigma$. By the compactness of $S$, we have
$$\Hom{}{S}{B^{(I)}}\cong \coprod_{I}\Hom{}{S}{B}$$
and $B^{(I)}$ is generated by $|I|$ copies of $\sigma$. It follows that $\Hom{}{S}{B^{(I)}}$ is generated by the family $\{\sigma\circ q_{i}\}_{i\in I}$. Therefore, we now only need to check that $\theta$ and the map, $\Hom{}{S}{\langle\phi (q_{i})\rangle}$, induced by $\langle\phi(q_{i})\rangle$ coincide on this set of generators.

By the commutativity of diagram \eqref{comm_square} we have:
\begin{eqnarray*}
\theta(q_{i}\circ\sigma) & = & \Hom{}{\sigma}{B^{(J)}}\circ\phi(q_{i}) \\
			 & = & \phi(q_{i})\circ\sigma \\
			 & = & (\langle\phi(q_{i})\rangle\circ q_{i})\circ\sigma \\
			 & = & \langle\phi(q_{i})\rangle\circ(q_{i}\circ\sigma) \\
			 & = & \Hom{}{S}{\langle\phi(q_{i})\rangle}(q_{i}\circ\sigma)
\end{eqnarray*}
Hence, $\theta$ and $\Hom{}{S}{\langle\phi(q_{i})\rangle}$ coincide on a basis of $\Hom{}{S}{B^{(I)}}$, thus
$$\theta = \Hom{}{S}{\langle\phi(q_{i})\rangle}$$
with $\langle\phi(q_{i})\rangle\in\Hom{}{B^{(I)}}{B^{(J)}}$. Therefore, the functor $\Hom{}{S}{-}$ is full and faithful.
\epf

\thm
Under the hypotheses of Setup \ref{setup}, the functor $$\Hom{}{S}{-}:\mathcal{C}\rightarrow \textnormal{Mod}(k^\textnormal{op}),$$ is an equivalence of categories, and hence, the coheart $\mathcal{C}$ of the non-degenerate co-$t$-structure obtained in Theorem \ref{prop4} is an abelian category.
\label{theorem}
\ethm

\pf
By Lemma \ref{surjects} and Proposition \ref{fullyfaithful}, $\Hom{}{S}{-}$ is dense and fully faithful. Hence, by \cite[Theorem IV.4.1]{MacLane}, $\Hom{}{S}{-}$ is an equivalence of categories.
\epf

\

Theorems \ref{prop4} and \ref{theorem} now combine to give Theorem \ref{cochaintheorem}.

Although it is known that the coheart of a co-$t$-structure is not always an abelian subcategory of $\T$, see \cite{Bondarko}, Theorem \ref{theorem} leads us to pose the following question.

\newtheorem{question}[dfn]{Question}
\begin{question}
Under what circumstances is the coheart of a co-$t$-structure on a triangulated category $\T$ an abelian subcategory of $\T$?
\end{question}

\noindent
\textbf{Acknowledgment.}  The author would like to thank his supervisor, Peter J\o rgensen, for all the help and advice he has given during the preparation of this paper, and also to thank the University of Leeds and EPSRC of the United Kingdom for financial support. In addition, the author is particularly grateful for the useful comments made by the referees.


\begin{thebibliography}{99}

\bibitem{Aldrich}
S. T. Aldrich and J. R. Garc\'ia Rozas, \textit{Exact and semisimple differential graded algebras}; Comm. Algebra \textbf{30} (2002), 1053-1075.

\bibitem{Avramov}
L. Avramov and S. Halperin, \textit{Through the looking glass: a dictionary between rational homotopy theory and local algebra};  Algebra, algebraic topology and their interactions (Stockholm, 1983),  1--27, Lecture Notes in Math., 1183, Springer, Berlin, 1986. 

\bibitem{Beilinson}
A. A. Beilinson, J. Bernstein and P. Deligne, \textit{Faisceaux pervers}; Ast\'erique \textbf{100} (1982).

\bibitem{Beligiannis_Reiten}
A. Beligiannis and I. Reiten, ``Homological and Homotopical Aspects of Torsion Theories''; preprint (2002), to appear in Mem. Amer. Math. Soc.

\bibitem{Bergh}
M. van den Bergh, \textit{A Remark on a Theorem by Deligne};  Proc. Amer. Math. Soc.  \textbf{132}  (2004),  no. 10, 2857--2858.

\bibitem{Bernstein}
J. Bernstein and V. Lunts, ``Equivariant Sheaves and Functors''; Lecture Notes in Mathematics 1578, Springer-Verlag, Berlin Heidelberg, 1994.

\bibitem{Bondarko}
M. V. Bondarko, \textit{Weight structures for triangulated categories: weight filtrations, weight spectral sequences and weight complexes; applications to motives and to the stable homotopy category}; arXiv:0704.4003v1 [math.KT] 30 Apr 2007.

\bibitem{Enochs}
E. Enochs, \textit{Injective and flat covers, envelopes and resolvents}; Israel J. Math. \textbf{39} (1981), 189-209.

\bibitem{Felix}
Y. F\'elix, S. Halperin and J-C. Thomas, ``Rational Homotopy Theory''; Graduate Texts in Mathematics 205, Springer, New York, 2001.

\bibitem{FJ2}
A. Frankild and P. J\o rgensen, \textit{Homological Identities for Differential Graded Algebras}; Journal of Algebra \textbf{265} (2003), pp. 114-135.

\bibitem{Miyachi}
M. Hoshino, Y. Kato and J. Miyachi, \textit{On $t$-structures and torsion theories induced by compact objects}; J. Pure Appl. Algebra \textbf{167} (2002), 15-35.

\bibitem{Iyama}
O. Iyama and Y. Yoshino, \textit{Mutations in triangulated categories and rigid Cohen-Macaulay modules}; arXiv:math.RT/0607736 v1 28 Jul 2006.

\bibitem{Kashiwara}
M. Kashiwara and P. Schapira, ``Sheaves on Manifold''; A Series of Comprehensive Studies in Mathematics 292, Springer-Verlag, Berlin Heidelberg, 1990.

\bibitem{MacLane}
S. MacLane, ``Categories for the Working Mathematician''; Springer, New York, 1971.

\bibitem{Neeman}
A. Neeman, \textit{The Grothendieck duality theorem via Bousfield's techniques and Brown representability}; J. Amer. Math. Soc. \textbf{9} (1996),  205-236.

\bibitem{Neeman2}
A. Neeman, ``Triangulated Categories''; Annals of Mathematics Studies, Princeton University Press, Princeton and Oxford, 2001.

\bibitem{Rotman}
J. J. Rotman, ``An Introduction to Algebraic Topology''; Graduate Texts in Mathematics 119, Springer-Verlag, New York, 1988.

\end{thebibliography}
\end{document}